\documentclass{gtart_h}

\def\ifplaintex{\expandafter\ifx\csname documentclass\endcsname\relax}

\def\gtp{{\mathsurround=0pt\it $\cal G\mskip-2mu$eometry \&\ 
$\cal T\!\!$opology $\cal P\!$ublications}}  

\def\recd{{\small Received:\qua\receiveddate\ifx\reviseddate\relax
\else\qquad Revised:\qua\reviseddate\fi\par}} 

\def\lognumber#1{\def\thelognumber{#1}}
\def\volumenumber#1{\def\thevolumenumber{#1}}
\def\volumeyear#1{\def\thevolumeyear{#1}}
\def\papernumber#1{\def\thepapernumber{#1}}
\def\pagenumbers#1#2{\def\startpage{#1}\def\finishpage{#2}}
\def\published#1{\def\publishdate{#1}}

\def\received#1{\def\receiveddate{#1}}

\def\accepted#1{\def\accepteddate{#1}}
\def\asciititle#1{\def\theasciititle{#1}}

\long\def\asciiabstract#1{\long\def\theasciiabstract{#1}}
\def\asciikeywords#1{\def\theasciikeywords{#1}}

\let\\\par\let\thelognumber\relax\let\thevolumenumber\relax
\let\thepapernumber\relax\let\thevolumeyear\relax\let\startpage\relax
\let\finishpage\relax\let\publishdate\relax\let\receiveddate\relax
\let\reviseddate\relax\let\accepteddate\relax\let\theasciititle\relax
\let\theasciiauthors\relax
\let\theasciiabstract\relax\let\theasciikeywords\relax

\let\theasciiemail\relax

\ifplaintex
\font\logobig=cmssbx10 scaled 3836
\font\logomed=cmssbx10 scaled 2557
\else
\font\logobig=cmssbx10 scaled 4200
\font\logomed=cmssbx10 scaled 2800
\fi

\long\def\makeagttitle{   
\count0=\startpage
\agt\hfill      
\hbox to 45truept{\vbox to 0pt{\vglue -13truept{\logomed A\kern -.37em{\logobig 
T}\kern -.38em G}\vss}\hss}
\break
{\small Volume \thevolumenumber\ (\thevolumeyear)
\startpage--\finishpage\nl
Published: \publishdate}

\vglue .25truein

{\parskip=0pt\leftskip 0pt plus
1fil\def\\{\par\smallskip}{\Large\bf\thetitle}\par\medskip} \vglue
0.05truein

{\parskip=0pt\leftskip 0pt plus 1fil\def\\{\par}{\sc\theauthors}
\par\medskip}%
 
\vglue 0.03truein 

{\small\leftskip 25truept\rightskip 25truept{\bf Abstract}\stdspace\theabstract

{\bf AMS Classification}\stdspace\theprimaryclass
\ifx\thesecondaryclass\relax\else; \thesecondaryclass\fi\par
{\bf Keywords}\stdspace \thekeywords\par}\vglue 7truept

}   

\ifplaintex
\font\phead=cmsl9 scaled 950
\font\pnum=cmbx10 scaled 913
\font\pfoot=cmsl9 scaled 950
\headline{\vbox to 0pt{\vskip -4.5mm\line{\small\phead\ifnum
\count0=\startpage ISSN 1472-2739 (on-line) 1472-2747 (printed)
\hfill {\pnum\folio}\else\ifodd\count0\def\\{ }%
\ifx\theshorttitle\relax\thetitle\else\theshorttitle\fi\hfill{\pnum\folio}
\else\def\\{ and }{\pnum\folio}\hfill\ifx\theshortauthors\relax\theauthors
\else\theshortauthors\fi\fi\fi}\vss}}
\footline{\vbox to 0pt{\vglue 0mm\line{\small\pfoot\ifnum\count0=\startpage
\copyright\ \gtp\hfill\else
\agt, Volume \thevolumenumber\ (\thevolumeyear)\hfill\fi}\vss}}
\else
\font=cmsl9 scaled 1050
\font\lnum=cmbx10 
\font=cmsl9 scaled 1050
\makeatletter
\def\@oddhead{{\small\lhead\ifnum\count0=\startpage ISSN 1472-2739 
(on-line) 1472-2747 (printed)\hfill {\lnum\number\count0}\else\ifodd\count0
\def\\{ }\ifx\theshorttitle\relax \thetitle \else\theshorttitle\fi\hfill
{\lnum\number\count0}\else\def\\{ and }{\lnum\number\count0}
\hfill\ifx\theshortauthors\relax 
\theauthors\else\theshortauthors\fi\fi\fi}}\def\@evenhead{\@oddhead}
\def\@oddfoot{\small\lfoot\ifnum\count0=\startpage\copyright\ \gtp\hfill\else
\agt, Volume \thevolumenumber\ (\thevolumeyear)\hfill\fi}
\def\@evenfoot{\@oddfoot}
\makeatother
\fi
\let\maketitlepage\makeagttitle
\let\makeshorttitle\maketitlepage
\let\maketitle\maketitlepage

\newwrite\gtoutfile
\long\gdef\makeheadfile{  
{\def\\{, }\def\s{ }
\immediate\openout\gtoutfile head.xxx
\immediate\write\gtoutfile{Proxy-for: \ifx\theasciiauthors\relax
\theauthors\else\theasciiauthors\fi\s<\ifx\theasciiemail\relax\theemail\else\theasciiemail\fi>}
\immediate\write\gtoutfile{\noexpand\\}
\immediate\write\gtoutfile{Authors: \ifx\theasciiauthors\relax
\theauthors\else\theasciiauthors\fi}
{\def\\{ }\immediate\write\gtoutfile{Title: \ifx\theasciititle\relax
\thetitle\else\theasciititle\fi}}
\immediate\write\gtoutfile{Subj-class: GT or SG, GR etc}
\immediate\write\gtoutfile{MSC-class: \theprimaryclass\ifx\thesecondaryclass\relax\else, \thesecondaryclass\fi}
\immediate\write\gtoutfile{Journal-ref: Algebr. Geom. Topol. \thevolumenumber\s
(\thevolumeyear) \startpage-\finishpage}
\immediate\write\gtoutfile{Comments: Published by Algebraic and
Geometric Topology at}
\immediate\write\gtoutfile{\s\s\s  http://www.maths.warwick.ac.uk/agt/AGTVol\thevolumenumber/agt-\thevolumenumber-\thepapernumber.abs.html}
\immediate\write\gtoutfile{\noexpand\\}
\immediate\write\gtoutfile{}
\ifx\theasciiabstract\relax
\immediate\write\gtoutfile{\theabstract}\else
\immediate\write\gtoutfile{\theasciiabstract}\fi
\immediate\write\gtoutfile{}
\immediate\write\gtoutfile{\noexpand\\}
\immediate\write\gtoutfile{}
\immediate\closeout\gtoutfile}}  

\def\maketitlepage{\makeagttitle\makeheadfile}
\let\makeshorttitle\maketitlepage
\let\maketitle\maketitlepage

\lognumber{11}
\volumenumber{5}
\volumeyear{2005}
\papernumber{11}
\pagenumbers{207}{217}
\received{23 February 2005} 
\accepted{6 March 2005}
\published{28 March 2005}

\usepackage{amsmath,amssymb,graphicx}

\theoremstyle{plain}
\newtheorem{thm}{Theorem}[section]
\newtheorem{lem}[thm]{Lemma}
\newtheorem{cor}[thm]{Corollary}

\newcommand{\R}{\ensuremath{\mathbb{R}}}
\newcommand{\C}{\ensuremath{\mathbb{C}}}
\newcommand{\s}{\ensuremath{\mathrm{S}}}
\newcommand{\Z}{\ensuremath{\mathbb{Z}}}
\newcommand{\Q}{\ensuremath{\mathbb{Q}}}
\newcommand{\slc}{\mathrm{SL}_2(\C)}

\renewcommand{\L}{\mathcal{L}}

\newcommand{\inv}[1]{{#1}^{-1}}
\newcommand{\pie}{\pi_1}
\newcommand{\bound}{\partial}
\newcommand{\tr}{\mathrm{Tr}}
\newcommand{\AND}{\quad \mathrm{and} \quad}

\newcommand{\cc}{\mathcal{C}}
\renewcommand{\th}[1]{#1^ \mathrm{th}}
\def\co{\colon\thinspace}

\begin{document}

\title{All roots of unity are detected\\by the A--polynomial}
\asciititle{All roots of unity are detected\\by the A-polynomial}
\author{Eric Chesebro}
\address{Department of Mathematics, The University of Texas at 
Austin\\Austin, TX 78712-0257, USA}
\email{chesebro@math.utexas.edu}
\urladdr{http://www.ma.utexas.edu/users/chesebro/}
\begin{abstract}
For an arbitrary positive integer $n$, we construct infinitely many one-cusped hyperbolic 3--manifolds where each manifold's A--polynomial detects every $\th{n}$ root of unity.  This answers a question of Cooper, Culler, Gillet, Long, and Shalen as to which roots of unity arise in this manner.
\end{abstract}

\asciiabstract{%
For an arbitrary positive integer n, we construct infinitely many
one-cusped hyperbolic 3-manifolds where each manifold's A-polynomial
detects every n-th root of unity.  This answers a question of
Cooper, Culler, Gillet, Long, and Shalen as to which roots of unity
arise in this manner.}

\asciikeywords{Character variety, ideal point, A-polynomial}
\maketitle

\section{Introduction}

Let $N$ be a compact orientable 3--manifold whose boundary is a disjoint union of incompressible tori.  The seminal paper \cite{CS} shows that ideal points of the $\slc$ character variety give rise to non-trivial splittings of $\pie(N)$ and hence to essential surfaces in $N$.   If such a surface has non-empty boundary, we say that the ideal point \textit{detects a boundary slope}.  (See Subsection \ref{sec:char} for a more detailed discussion.)

It is shown in \cite{CCGLS}, that canonically associated to a boundary slope detected by an ideal point is a root of unity.  For if $\alpha$ represents the detected boundary slope, then the limiting value of the characters evaluated at $\alpha$ is of the form $\xi + \inv{\xi}$,  where $\xi$ is a root of unity.  We say that $\xi$ \textit{appears} at this ideal point.  It is a natural question to ask which roots of unity arise in this manner. 

Roots of unity that appear at ideal points carry topological information about associated surfaces.  Let $x$ be an ideal point which detects a boundary slope, $\xi$ be the associated root of unity, and $\Sigma$ be a surface associated to $x$.  Assume that the number of boundary components of $\Sigma$ is minimal over all such surfaces.  By Corollary 5.7 of \cite{CCGLS}, if $S$ is a component of $\Sigma$ then  the order of $\xi$ divides the number of boundary components of $S$. 

At present, little seems known about which roots of unity can occur.  If $N$ is the complement of a 2--bridge knot, then the essential surfaces in $N$ are classified in \cite{HT}.  In particular, any connected essential surface in $N$ has either one or two boundary components.  Therefore, the only roots of unity that appear in the character variety for $N$ are $\pm1$.  In fact, in most known examples, the roots of unity that appear are $\pm1$.  It is now known, however,  that higher order roots do appear.  In \cite{CLroots}, Cooper and Long mention an example with $\th{11}$ roots, Dunfield gives examples with $\th{4}$ and $\th{6}$ roots in \cite{Dnroots}, and Kuppum and Zhang give other examples with $\th{4}$ roots in \cite{ZK}.

We prove the following theorem.

\begin{thm} \label{th:first}
Let  $n$ be a positive integer.  Then there exist infinitely many one-cusped finite volume hyperbolic 3--manifolds $\{M_i\}_{i=1}^\infty$ where if $\xi$ is any $\th{n}$ root of unity then $\xi$ appears as a root of unity for every $M_i$. 
\end{thm}

The idea of the proof is this.  For each $n \in \Z^+$, we construct a 3--manifold $N_n$ for which $\th{n}$ roots of unity appear.  The manifolds $N_n$ have non-trivial JSJ decompositions; hence they are not hyperbolic.  We may, however, use standard techniques to construct hyperbolic 3--manifolds admitting degree-one maps onto $N_n$.  The root of unity property is inherited by the hyperbolic manifolds through the induced epimorphism on fundamental groups.

\section{Background and notation}

For more detailed discussions of topics in this section, see \cite{CCGLS} and \cite{CS}.

\subsection{Notation}
Throughout this paper all 3--manifolds are assumed to be connected and orientable.  We assume further that all boundary components of 3--manifolds are homeomorphic to tori.  All surfaces are assumed to be properly embedded and orientable.  

A \textit{slope} on a 3--manifold $N$ is the isotopy class of an unoriented essential simple closed curve on $\bound N$.  Every slope corresponds to a pair $\{ \alpha, \inv{\alpha} \}$ of primitive elements in $\pie(\bound N)$.  A \textit{boundary slope} is a slope that is represented by a boundary component of an essential surface in $N$.  

Let $L$ be a disjoint union of smooth, simple, closed curves in $\s^3$.  We define the \textit{exterior of $L$} as the compact 3--manifold $E(L)=\s^3-\eta(L)$ where $\eta(L)$ is an open tubular neighborhood of $L$ in $\s^3$.  

\subsection{Character varieties and detected surfaces} \label{sec:char}

Let $\Gamma$ be a finitely generated group.  We write $X(\Gamma)$ to denote the set of characters of $\slc$ representations.  $X(\Gamma)$ is naturally a closed algebraic set in $\C^m$.  If $X$ is an irreducible algebraic curve in $X(\Gamma)$ then there is a smooth projective curve $\widetilde{X}$ and a birational map $f \co \widetilde{X}\longrightarrow X$.  The points $\widetilde{X}-\inv{f}(X)$ are called \textit{ideal points}.  Associated to each ideal point is a valuation which in turn gives rise to a non-trivial simplicial action on a tree.  If $\Gamma$ is the fundamental group of a 3--manifold $N$, we write $X(N)=X\bigl(\pie(N)\bigr)$.  In this case, the action associated to a ideal point may be used to construct essential surfaces in $N$.  We say that such a surface is \textit{associated to the action given by the ideal point.}

For any $\gamma \in \pie(\bound N)$,  there is a regular function $I_\gamma \co X \longrightarrow \C$ given by $I_\gamma(\chi)=\chi(\gamma)$.  These functions induce meromorphic functions on $\widetilde{X}$.  If $x$ is an ideal point, we consider the value $I_\gamma(x) \in \C \cup \{ \infty \}$.   If there exists a $\gamma \in \pie( \bound N )$ so that $I_\gamma$ has a pole at $x$ then there is a unique boundary slope $s=\{ \alpha, \alpha^{-1}\}$ with the property that $I_{\alpha^{\pm1}}(x)$ is finite.  Any essential surface $\Sigma$ coming from this ideal point will have non-empty boundary and any boundary component of $\Sigma$ represents the boundary slope $s$. 

It is shown in \cite{CCGLS} and \cite{CLroots} that the bounded value $I_{\alpha^{\pm 1}}(x)$ is of the form $\xi + \inv{\xi}$, where $\xi$ is a root of unity. 

\subsection{The A--polynomial} \label{sec:Apoly}

In this section we define the A--polynomial $A_N(L,M)$ and sketch its connection to the root of unity phenomena described above.  The primary reference is \cite{CCGLS}.  

Throughout this section, we fix a generating set $\{ \lambda, \mu \}$ for $\pie( \bound N )$.  If $\rho$ is a $\slc$-representation of $\pie(\bound N)$ with $\tr\bigl(\rho(\mu)\bigr)\neq \pm 2$ then $\rho$ is conjugate to a diagonal representation.  We are concerned with characters of representations, so we assume that $\rho$ is diagonal.  Let
$$\rho(\mu) \ = \ \begin{pmatrix} M & 0 \\ 0 & 1/M \end{pmatrix} \AND \rho(\lambda) \ = \ \begin{pmatrix} L & 0 \\ 0 & 1/L \end{pmatrix}.$$
It is not difficult to see that $X(\bound N)$ is the quotient of $\C^\times \times \C^\times$ by the involution $(L,M) \mapsto (\inv{L}, \inv{M})$. 

The inclusion $i \co \bound N \hookrightarrow N$ induces a regular map $i^\ast \co X(N) \longrightarrow X(\bound N)$ between the character varieties.  Let $\{X_j \}_{j=1}^k$ be the set of irreducible components of $X(N)$ with the property that $\dim_\C \bigl(i^\ast(X_j)\bigr) =1$ and every $X_j$ contains the character of an irreducible representation.  Define
$$ D_N \ = \ \bigcup_{j=1}^k i^\ast (X_j)$$
and $\overline{D_N}$ as the inverse image of $D_N$ under the quotient map $\C^\times \times \C^\times \longrightarrow X(\bound N)$.   If $A_N(L,M)\in \Q[L,M]$ is a defining polynomial for $\overline{D_N}$, it can be normalized to have relatively prime $\Z$--coefficients.  This is the A--polynomial for $N$.

The \textit{Newton polygon} for $A_N(L,M)= \sum a_{ij}L^i M^j$ is the convex hull of the set $\{ (i, j) \in \R^2 \, | \, a_{ij} \neq 0 \}.$  If $e$ is an edge of the Newton polygon for $A_N(L,M)$ and $p/q$ is its slope, we can change our choice of generating set for $\pie(\bound N)$ to $\{ \mu^p \lambda^q,  \mu^a \lambda^b \}$, where $a$ and $b$ are chosen so that $\det \bigl( \begin{smallmatrix} p & a \\ q & b \end{smallmatrix} \bigr) = 1.$  We now compute the A--polynomial which corresponds to the new generating set by setting
$$M \ = \ Z^bY^{-q} \AND L \ = \ Z^{-a}Y^p$$
in $A_N(L,M)$ and then normalizing as before.  The result is of the form
$$A_N(Y,Z) \ = \ f(Z) + E(Y,Z).$$
where $f$ is not identically zero, has at least one non-zero root $\xi$, and $E(0,Z)=0$ for every $Z \in \C$.  The \textit{edge polynomial} for $e$ is defined as 
$$f_e(t) \ = \ t^{-m} \cdot f(t),$$
where $m \in \Z^+$ is chosen so that $f_e \in \Z [t]$ and $f_e(0) \neq 0$.

Notice that the solution $(0,\xi)$ to $A_N(Y,Z)=0$ does not correspond to the image of a character $\chi_\rho$.  It does, however, correspond to the limit of a sequence of characters of irreducible representations $\bigl\{\chi_i\bigr\}_{i=1}^\infty$.  Furthermore, $I_{\mu^a \lambda^b}$ has a pole at this ideal point and the values $\bigl\{ I_{\mu^p \lambda^q}(\chi_i) \bigr\}$ converge to $\xi + \inv{\xi}$.  As outlined in Section \ref{sec:char}, $p/q$ is a boundary slope with respect to the basis $\{ \lambda, \mu \}$ and $\xi$ is a root of unity.

\section{Construction of the manifolds $N_n$}

Let $K$ be the figure-eight knot, $N_1=E(K)$, and $\Gamma_1=\pie(N_1)$.
We have the following presentation,
$$\Gamma_1 \  = \ \big\langle \, \mu,y \, \big| \, \mu W=Wy \, \big\rangle,$$
where $W=\inv{y}\mu y\inv{\mu}$.  A meridian for $N_1$ is $\mu$ and the longitude, $\lambda$, is given by $WW^\ast$ where $W^\ast=\inv{\mu} y \mu \inv{y}$. 

The A--polynomial for $N_1$, with respect to the basis $\{ \lambda, \mu \}$ for $\pie( \bound N_1)$ is 
\begin{equation} \label{eq:Apoly}
A_{N_1}(L,M) \ = \ -L + LM^2+M^4+2LM^4+L^2M^4+LM^6-LM^8
\end{equation}
(see \cite{CCGLS}).  The Newton polygon for the above polynomial has an edge $e$ with slope $4/1$.  As it will be needed in what follows, we note that the root of unity associated to this slope is $+1$.  To see this, set
$$M \ = \ ZY^{-1} \AND L \ = \ Z^{-3}Y^4$$
in \eqref{eq:Apoly}.  Multiply by $Z^3Y^4$ and normalize signs to get
$$ A_{N_1}(Y, Z) \ = \ (Z^7 - Z^8) + Y^2 (-Y^6 + Y^4Z^2 + 2Y^2Z^4 +Y^6Z +Z^6).$$
Therefore, the edge polynomial is 
$$f_e(t) \ = \ t^{-7}( t^7 - t^8) \ = \ 1-t, $$
as required.

\begin{figure}[ht!]
\begin{center}
\begin{picture}(0,0)%
\includegraphics[scale=.3333]{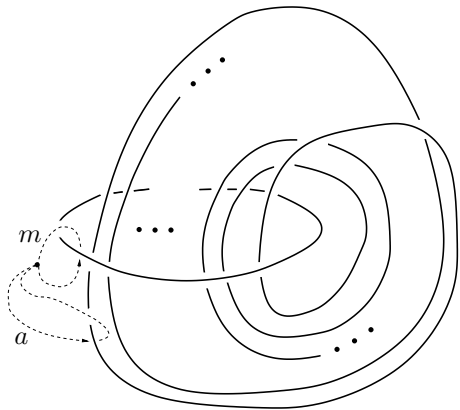}%
\end{picture}%
\setlength{\unitlength}{1316sp}
%
\begin{picture}(8619,7630)(115,-6020)
\put(205,-2855){\makebox(0,0)[lb]{\small$m$}}
\put(131,-4788){\makebox(0,0)[lb]{\small$a$}}
\end{picture}%
\end{center}
\caption{The link $L_n$} \label{Ln}
\end{figure}

\begin{figure}[ht!]
\begin{center}
\begin{picture}(0,0)%
\includegraphics[scale=.3333]{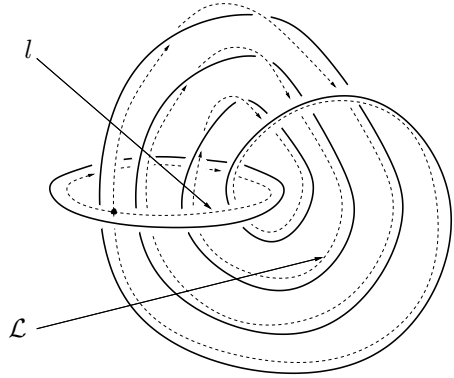}%
\end{picture}%
\setlength{\unitlength}{1316sp}
%
\begin{picture}(8610,7017)(436,-6333)
\put(751,-361){\makebox(0,0)[lb]{\small$l$}}
\put(451,-5611){\makebox(0,0)[lb]{\small$\mathcal{L}$}}
\end{picture}%
\caption{The link $L_4$ labeled with representatives for $l$ and $\L$}\label{L4}
\end{center}
\end{figure}

Now, fix an $n \in \Z^+$ and let $N_2^n=E(L_n)$, where $L_n$ is the $(2,2n)$--torus link shown in Figure \ref{Ln}.  Let $\Gamma_2^n=\pie(N_2^n)$.  Then,
$$\Gamma_2^n \ =  \ \big\langle \, a,m \, \big| \, a(am)^n=(am)^n a \, \big\rangle.$$
Let $\L \in \Gamma_2^n$ be the peripheral element which is homologous to $m^n$ and let $l \in \Gamma_2$ be the peripheral element homologous to $a^n$. See Figure \ref{L4} for the case when $n=4$.

We form the 3--manifold $N_n$ as the union of $N_1$ and $N_2^n$ by identifying the boundary of $N_1$ with one component of the boundary of $N_2$ by insisting that $l$ is identified with $\mu^4\lambda$ and $m$ is identified with $\mu^3\lambda$.  Then $\Gamma_n=\pie(N_n)$ is a free product with amalgamation $\Gamma_1 \ast _H \Gamma_2^n$ where $H = \langle  \lambda, \mu \rangle = \langle l, m \rangle$.   Furthermore, $\{\L, a \}$ is a generating set for $\pie(\bound N_n)$.

 \noindent \textbf{Remark}\qua  We chose the figure-eight knot and the links $L_n$ for simplicity, but there are many other choices that will also work.  For example, in defining $N_1$, we can choose any knot $K$ which has a detected boundary slope $s$ with associated root of unity $+1$.  In this case, we form the manifolds $N_n$ by identifying the curve $l$ to the boundary slope $s$.  Similarly, in choosing the manifolds $N_2^n$, we have considerable freedom.  We need only require that the curve which is identified to the detected boundary slope $s$ on $N_1$ is homologous to to the $\th{n}$ power of a primitive non-peripheral element of $\pie(\bound N_n)$.

\section{$N_n$ has $\th{n}$ roots}

\begin{thm} \label{th:main}
Fix $n \in \Z^+$ and let $\xi$ be any $n^\mathrm{th}$ root of unity.  Then there exists a sequence of irreducible representations 
$$ \bigl\{ \rho_i \co \Gamma_n \longrightarrow \slc \bigr\}_{i=1}^\infty$$
such that
$$\rho_i(\L) \ = \ \begin{pmatrix} B_i & 0 \\ 0 & B_i^{-1} \end{pmatrix} \AND \rho_i(a) \ = \ \begin{pmatrix} A_i & 0 \\ 0 & A_i^{-1} \end{pmatrix}.$$
Where, as $i \rightarrow \infty$, $B_i \rightarrow 0$ and $A_i \rightarrow \xi$.  In particular, $\xi$ appears at an ideal point for $N_n$.
\end{thm}

\textbf{Proof}\qua  To prove the theorem, we exhibit a sequence of representations which are irreducible when restricted to $\Gamma_1$, abelian when restricted to $\Gamma_2^n$, and that satisfy the conclusions of the theorem.

Recall the A--polynomial $A_{N_1}(L,M)$ given in \eqref{eq:Apoly} of Section \ref{sec:Apoly}.  Then for all but finitely many solutions $(L,M)$ to the equation $A_{N_1}(L,M) = 0$, we have an irreducible representation $\psi_1 \co \Gamma_1 \longrightarrow \slc$ with
$$\lambda \mapsto \begin{pmatrix} L & 0 \\ 0 & 1/L \end{pmatrix} \AND \mu \mapsto \begin{pmatrix} M & 0 \\ 0 & 1/M \end{pmatrix}.$$
The abelianization of $\Gamma_2^n$ is isomorphic to $\Z \times \Z$, so we get an abelian representation for $\Gamma_2^n$ by chosing any two diagonal elements of $\slc$.   More precisely, if $p, A \in \C^\times$ then there exists an abelian representation $\psi_2 \co \Gamma_2^n \longrightarrow \slc$ given by 
 $$m \mapsto  \begin{pmatrix} p & 0 \\ 0 & 1/p \end{pmatrix} \AND a \mapsto  \begin{pmatrix} A & 0 \\ 0 & 1/A \end{pmatrix}.$$  
Since $\psi_2$ is abelian and $\L$ and $l$ are homologous to $m^n$ and $a^n$ respectively, we have
$$\psi_2(\L) \ = \ \begin{pmatrix} p^n & 0 \\ 0 & p^{-n} \end{pmatrix} \quad \mathrm{and} \quad \psi_2(l) \ = \ \begin{pmatrix} A^n & 0 \\ 0 & A^{-n} \end{pmatrix}.$$
The mapping $\rho$  given by
$$ 
 \begin{array}{lcl}
 \rho(\mu) \ = \ \psi_1(\mu) & \qquad &
 \rho(y) \ = \ \psi_1(y) \\ 
 \rho(a) \ = \ \psi_2(a) & &
 \rho(m) \ = \ \psi_2(m)
 \end{array}
$$
 extends to an irreducible representation $\rho \co \Gamma_n \longrightarrow \slc$  if and only if
$$ \psi_1(\mu^4\lambda) \ = \ \psi_2(l) \quad \mathrm{and} \quad \psi_1(\mu^3 \lambda) \ = \ \psi_2(m).$$
Note that $\mu = l \inv{m}$ and $\lambda=l^{-3}m^4$ in $\Gamma$. Therefore, $\rho$ is a representation if and only if
$$ M \ = \ \frac{A^n}{p} \quad \mathrm{and} \quad L \ = \ \frac{p^4}{A^{3n}}.$$
If we replace the indeterminants $L$ and $M$ in the equation $A_{N_1}(L,M)=0$  with the above expressions, then after simplifying and multiplying by $A^{3n}p^4$ we get
\begin{equation} \label{eq:solns}
-p^8 + p^6 A^{2n} + A^{7n} + 2 p^4 A^{4n} + p^8 A^n + p^2 A^{6n} - A^{8n} \ = \ 0.
 \end{equation}
By construction, all but finitely many solutions to \eqref{eq:solns} correspond to representations of $\Gamma_n$.  They are irreducible because their restrictions to $\Gamma_1$ are irreducible, and so we have an irreducible representation $\rho \co \Gamma \longrightarrow \slc$ with 
$$a \mapsto \begin{pmatrix} A & 0 \\ 0 & 1/A \end{pmatrix} \AND \L \mapsto \begin{pmatrix} p^n & 0 \\ 0 & p^{-n} \end{pmatrix} $$
for all but finitely many solutions to \eqref{eq:solns}.

Let $\xi$ be any $n^\mathrm{th}$ root of unity.  Note that setting $p=0$ and $A=\xi$ we get a solution to \eqref{eq:solns}.  This solution does not correspond to a representation of $\Gamma$, but we can take a sequence of distinct solutions  $\bigl\{ (p_i, A_i)\bigr\}_{i=1}^\infty$ which do correspond to irreducible representations and converge to $(0,\xi)$. Let $B_i =p_i^n$.  

This gives us a sequence $\{\rho_i \}_{i=1}^\infty$ of representations that satisfy the conclusions of the theorem.  Note that $\tr \bigl( \rho_i(a) \bigr)$ tends toward $\xi + \inv{\xi}$ and $\Big| \tr \bigl( \rho_i(\L) \bigr) \Big|$ tends to $\infty$.  Therefore, the boundary slope $\{ a, \inv{a} \}$ is detected and $\xi$ appears at the associated ideal point.\endproof

The root of unity behavior demonstrated above is reflected in the A--polynomial for $N_n$.
 
\begin{cor}
If $A_{N_n}(X,Y)$ is the A--polynomial for $N_n$, with respect to the generating set $\{ \L, a \}$ for $\pie(\bound N_n)$, then the Newton polygon for $A_{N_n}$ has a vertical edge $e$ such that $1-t^n$ divides the associated edge polynomial $f_e(t).$
\end{cor}

\textbf{Proof}\qua Given any $\th{n}$ root of unity, let $\{\rho_i \}_{i=1}^\infty$ be a sequence of representations as given by Theorem \ref{th:main}.  Recall that $\Big| \tr \bigl( \rho_i(\L) \bigr) \Big|$ tends to $\infty$.  Therefore, the set $\big\{ i^\ast (\rho_i) \big\}_{i=1}^\infty \subset D_{N_n}$ contains an infinite number of distinct points.  Then by passing to a subsequence we may assume $\big\{ i^\ast (\rho_i) \big\}_{i=1}^\infty$ is a set of distinct points on a single curve component $\cc$ of $D_{N_n}$.  Recall that the choice of generating set $\{\L, a \}$ induces a quotient $\C^\times \times \C^\times \longrightarrow X(\bound N_n)$.  Let $\overline{\cc} \subset \overline{D_{N_n}} \subset \C^\times \times \C^\times$ be the inverse image of $\cc$ under this quotient map.  Let $P(X,Y)$ be the irreducible polynomial that defines the curve $\cc$.  

By construction, $P(0,\xi)=0$, however $P(0,Y)$ is not identically zero since $P(X,Y)$ is irreducible.  Therefore, $P(0,Y)$ is a non-trivial, single variable polynomial.  Also, by definition of the A--polynomial, $P(X,Y)$ divides $A_{N_n}(X,Y)$.  So $P(0,Y)$ divides $A_{N_n}(0,Y)$.  This implies that the Newton polygon for $A_{N_n}$ has a vertical edge $e$ and $f_e(\xi)=0$.  Since this is true for any $\th{n}$ root of unity, we conclude that $1-t^n$ divides $f_e$.  \endproof

\section{Hyperbolic examples} \label{sec:hyp}

Here we extend the result to hyperbolic examples and complete the proof of Theorem \ref{th:first}.

Before we prove the theorem, we must state two lemmas.  First, as a special case of the main theorem in \cite{Mexcel}, we have
\begin{lem} \label{lem:one}
For each $n \in \Z^+$, $N_n$ contains infinitely many null-homotopic hyperbolic knots.
\end{lem}

Our second lemma is standard, for example, see Proposition 3.2 in \cite{BW}.  The argument given there is for closed manifolds, but it follows verbatim in our setting.

\begin{lem} \label{lem:two}
Let $k$ be a null homotopic knot in $N_n$ and $M$ be a 3--manifold obtained by a Dehn surgery on $k$.  Then there is a degree-one map $f \co M \longrightarrow N_n$.
\end{lem}

\textbf{Proof of Theorem \ref{th:first}}\qua   Let $k$ be a knot in $N_n$ given by Lemma \ref{lem:one}.  By Thurston's hyperbolic Dehn surgery theorem \cite{Th}, for all but finitely many slopes $s$, the manifold obtained by doing Dehn surgery on $k$ with surgery slope $s$ is hyperbolic.  This gives us an infinite set of distinct hyperbolic manifolds $\{ M_i \}_{i=1}^\infty$.  

By Lemma \ref{lem:two}, for every $i \in \Z^+$, we have a degree-one map $f_i \co M_i \longrightarrow N_n$.  Since $f_i$ is degree-one, the induced map $f_i ^\ast \co \pie(M_i) \longrightarrow \Gamma_n$ is an epimorphism.  Furthermore, the restriction to the peripheral subgroup of $M_i$ is an epimorphism onto the peripheral subgroup of $N_n$. 

Let $\alpha$ be a primitive element of $\pie(\bound M_i)$ so that $f_i^\ast(\alpha) = a$.  Similarly, take $\beta$ to be a primitive element of $\pie(\bound M_i)$ so that  $f_i^\ast(\beta) = \L$.  Pick $\xi$ to be any $\th{n}$ root of unity.  Then Theorem \ref{th:main} gives us a sequence of irreducible representations $\big\{\rho_j \co \Gamma_n \longrightarrow \slc \big\}_{j=1}^\infty$ such that 
$$\rho_i(\L) \ = \ \begin{pmatrix} B_i & 0 \\ 0 & B_i^{-1} \end{pmatrix} \AND \rho_i(a) \ = \ \begin{pmatrix} A_i & 0 \\ 0 & A_i^{-1} \end{pmatrix}$$
where $B_i \rightarrow 0$ and $A_i \rightarrow \xi$.  Composing the representations with $f_i^\ast$, we get the sequence  $\{ \rho_j \circ f_i^\ast \}_{j=1}^\infty$ of irreducible representations for $\pie(M_i)$.  Furthermore,
$$\lim_{j \to \infty} \tr\big(  \rho_j \circ f_i^\ast (\alpha)\big) \ = \ \lim_{j \to \infty} \big(A_j + A_j^{-1}\big) \ = \ \xi + \inv{\xi}$$
and 
$$\lim_{j \to \infty} \Big|\tr\big(  \rho_j \circ f_i^\ast (\beta)\big)\Big| \ = \ \lim_{j \to \infty} \big|B_j + B_j^{-1}\big| \ = \ \infty.$$ 
Therefore, $\alpha$ represents a detected boundary slope and the root of unity that appears at the associated ideal point is $\xi$.  \endproof

\section{Computations for $n=5$}

We compute the factor of the A--polynomial for $N_5$ which corresponds to the representations that are irreducible on $N_1$ and abelian on $N_2$.  As in the proof of Theorem \ref{th:main}, we want to replace the variables $L$ and $M$ in \eqref{eq:Apoly} with the variables $A$ and $B$ which represent the eigenvalues of $\rho(a)$ and $\rho(\L)$ respectively. Here we take resultants instead of directly substituting expressions in terms of $p$ and $A$.  If we substitute directly, we will get an expression with non-integral exponents.  The following calculations are done with a computer.

Starting with $A(L,M)$ we take a resultant with respect to $M-A^5 \inv{p}$ to eliminate the indeterminant $M$.  Next we take a resultant with respect to $L-p^4 A^{-15}$ to eliminate $L$, and finally with respect to $B-p^5$ to eliminate $p$.  We conclude that the polynomial equation, 
\begin{multline} \nonumber
0 \ = \  -A^{175} + 5 A^{180} - 10 A^{185} + 10 A^{190} - 5 A^{195} + A^{200} + 5 A^{135} B^2 - 35 A^{140} B^2 + \\
\qquad \qquad + 60 A^{145} B^2 - 26 A^{150} B^2 - 5 A^{155} B^2 - 5 A^{90} B^4 -45 A^{95} B^4 +\\
\qquad \qquad + 98 A^{100} B^4 - 45 A^{105} B^4 - 5 A^{110} B^4 - 5 A^{45} B^6 -26 A^{50} B^6 +\\
\qquad \qquad + 60 A^{55} B^6 - 35 A^{60} B^6 + 5 A^{65} B^6 + B^8 - 5 A^5 B^8 + 10 A^{10} B^8 - \\
 - 10 A^{15} B^8 + 5 A^{20} B^8 - A^{25} B^8.
\end{multline}
holds true for all such representations.  Therefore the right hand side is a factor of the A--polynomial for $N_5$ with respect to the generating set $\{\L, a \}$.  Evaluating at $B=0$ we get,
\begin{eqnarray*} 
 0 & = & -A^{175} + 5 A^{180} - 10 A^{185} + 10 A^{190} - 5 A^{195} + A^{200} \\
&  = & A^{175}(A^5-1)^5.  
\end{eqnarray*}
Hence we conclude that there is a vertical edge $e$ in the Newton polygon for the A--polynomial of $N_5$.  Furthermore, if $f_e(t)$ is the associated edge polynomial then $(t^5-1)^5$ divides $f_e(t)$.  

\section{Final remarks}

It should be noted that although the manifolds constructed in Section \ref{sec:hyp} are hyperbolic, the ideal points at which $\th{n}$ roots of unity appear are not on components of the character varieties that correspond to the complete hyperbolic structure of the manifold.  This is because the components considered here arise from epimorphisms with non-trivial kernels and so the characters on these components are not characters of discrete faithful representations.  The question of whether all roots of unity appear at ideal points on hyperbolic components remains unanswered.

\section{Acknowledgments}

The author would like to thank his advisor, Alan Reid, for his encouragement and many helpful conversations.  He would also like to thank Richard Kent for pointing out reference \cite{Mexcel}.

\Addresses

\end{document}